\documentclass[12pt,oneside]{amsart}

\usepackage{amssymb}
\usepackage[latin1]{inputenc}

\def\N{\mathbb{N}}
\def\Z{\mathbb{Z}}

\def\e{\varepsilon}

\def\l{\lambda}

\def\w{\omega}
\def\G{\Gamma}

\def\cal{\mathcal}

\def\lgr{\mathrm{lgr}}

\newtheorem{thm}{Theorem}

\newtheorem{lem}{Lemma}

\newtheorem{pr}{Property}

\title{On stable norm in word hyperbolic groups}

\begin{document}
\maketitle 
\begin{center}
{\sc Jean-Philippe PR\' EAUX}\footnote[1]{Centre de recherche de l'Armée de l'Air, Ecole de l'air, F-13661 Salon de
Provence air}\ \footnote[2]{Centre de Math\'ematiques et d'informatique, Universit\'e de Provence, 39 rue
F.Joliot-Curie, F-13453 marseille
cedex 13\\
\indent {\it E-mail :} \ preaux@cmi.univ-mrs.fr\\
{\it Mathematical subject classification : 20F67}}
\end{center}
\begin{abstract}
This work is concerned with the stable norm in word hyperbolic groups, as defined  in \cite{gro}. We give a short
elementary proof of one of its basic property, that is existence of a computable non null uniform lower bound for
stable norm in a word hyperbolic group.
\end{abstract}

\section*{Introduction}
Though in word hyperbolic groups (cf. \cite{gro}, \cite{cdp}, \cite{gdlh} or \cite{sho}) the stable norm has appeared
to be a useful tool, some of its basic properties have only been sketched in papers like \cite{gro} and \cite{delz}, at
least as far as we're concerned. We are interested in the existence of a computable uniform lower bound for the stable
norm in a word hyperbolic group and give a short elementary complete proof.

\section{Definition of stable norm}

Given a  group $G$  on a finite generating set $S$, a word on $S\cup S^{-1}$ naturally represents an element of $G$. We
denote by $\lgr(.)$ its length as a word  and :
$$|g|=inf\{\lgr(\w)\,|\, \w\;\mathrm{is\; a\; word\; on}\ S\cup S^{-1}\ \mathrm{representing}\ g\}$$
which allows to define the word metric $d_S$ on $G$ by $d_S(g_1,g_2)=|g^{-1}g_2|$.\smallskip\\
\noindent{\bf Definition.} {\sl Let $G$ be a  group and $g\in G$ an element with infinite order, we denote by :}
$$[[g]]=\lim_{n\rightarrow +\infty} \frac{|g^n|}{n}$$
{\sl the {\bf stable norm} of the element $g$.}\smallskip\\
\indent The limit of $|g^n|/n$ exists : with the triangle inequalities, for any $n,p\in\N$, $|g^{n+p}| \leq
|g^n|+|g^p|$, that is the sequence $(|g^n|)_n$ is sub-additive.

\begin{lem}
Let $(u_n)_n$ be a sub-additive sequence   of positive real numbers, then the sequence $(u_n/n)_n$ converges.
\end{lem}

\noindent {\sl Proof.} Let $n\geq m > 0$ ; consider the euclidian division of $n$ by $m$ : $n=pm+r$ with $p,r\in\N$,
$0\leq r <m$. Since $(u_n)_n$ is sub-additive, $0\leq u_n\leq u_{pm}+u_r\leq p\,u_m+\max_{0\leq r <m} u_r$. Then :
$$0\leq\frac{u_n}{n}\leq \frac{pm}{n}\frac{u_m}{m}+\frac{1}{n}\max_{0\leq r <m} u_r
\leq \frac{u_m}{m}+\frac{1}{n}\max_{0\leq r <m} u_r$$ Making $n$ tend to $\infty$ one obtains for any $m\in\N^*$ :
$0\leq \lim\ \sup {u_n}/{n}\leq {u_m}/{m}$, so that $\lim\, \sup {u_n}/{n}$ is finite and moreover :
 $\lim\ \sup {u_n}/{n}\leq \inf_{m>0} {u_m}/{m}\leq \lim\
\inf {u_m}/{m}$, so that : $\lim\ \sup {u_n}/n\leq\lim\ \inf u_n/n$, which proves the assumption.\hfill$\square$

\section{Basic properties}
The two following are basic properties of stable norm :

\begin{pr}\label{1}
The stable norm is an invariant of conjugacy classes.
\end{pr}

\noindent {\sl Proof.} Suppose $u$ and  $v$ are infinite order elements of a  group $G$, and $u=ava^{-1}$ for some
$a\in G$. Then for any $n\in\N^*$, $u^n=av^na^{-1}$ which implies that $|\, |u^n|-|v^n|\, |\leq 2\,|a|$. Dividing by
$n$ and making $n$ tend to $\infty$ one obtains  $[[u]]=[[v]]$.\hfill$\square$\\

\begin{pr}\label{2}
For any $n\in \Z^*$, $[[g^n]]=|n|\, .[[g]]$.
\end{pr}

\noindent{\sl Proof.} Given $n\in\Z^*$,
$$[[g^n]]=\lim_{t\rightarrow +\infty} \frac{|g^{nt}|}{t}=\lim_{t\rightarrow +\infty}|n|\, \frac{|g^{|n|t}|}{|n|t}=
|n|\lim_{t\rightarrow +\infty}\frac{|g^{t}|}{t}=|n|\, .[[g]]$$
\hfill$\square$

\section{Uniform lower bound for stable norm}
We now give the proof of the property we are interested in.

\begin{thm}
Let $G$ be a word hyperbolic group ; there exists a computable constant $K>0$ such that for any infinite order element
$g$ of $G$, $[[g]]\geq K$.
\end{thm}

\noindent{\sl Proof.} Denote by $\G(G)$ the Cayley graph of $G$ respectively to a given finite generating set; $G$ is
given the word metric which makes it isometric to the vertex set of $\G(G)$.
A word on the generators of $G$ is said  {\sl cyclically reduced} if all its cyclic conjugates are geodesic words.
Obviously for any element $g\in G$ there exists a cyclically reduced word $\w$ representing an element in the conjugacy
class of $g$. Whenever $\w$ is a word representing an infinite order element of $G$, one can consider a segment in
$\G(G)$ with label $\w$ and the infinite path,  that we denote $[\w^-,\w^+]$, which is the orbit of this segment under
the action of the cyclic subgroup generated by $\w$.
If $\w$ is a cyclically reduced word of length $k>0$, then $[\w^-,\w^+]$ is a $k$-local geodesic, that is,  any of its
subpaths of length at most $k$ is a geodesic.

 Since the word problem is solvable in $G$
  one can algorithmically transform a given word $w$ into a cyclically reduced one lying in the same conjugacy class.
   We will  make use of this algorithm together with the fact (theorem 3.1.4, \cite{cdp}) that there
exists computable constant numbers $k>0$, $\l\geq 1$ and $\e\geq 0$ such that each $k$-local geodesic in the Cayley
graph of $G$ is a $(\l,\e)$-quasigeodesic. In the following we fix the constants $k$, $\l$, $\e$.

Denote by $\cal B$ the ball in $G$ with center $1$ and radius $k-1$. Apply the following process  for each element $g$
in $\cal B$ : Initially $i=2$ ; change $g^i$ into a cyclically reduced word in the same conjugacy class. If its length
is at least $k$ then stop with $n_g=i$ else if it is non null restart with $i+1$ instead of $i$, else stop.
 The process will terminate since either $g$ has finite order or
 $g^p$ cannot be conjugate with $g^q$ for $0<p<q$ (hyperbolic groups do not contain any Baumslag-Solitar group). Among all the integers $n_g$ obtained
for (infinite order) elements $g$, denote by  $n_{max}$ the supremum ; clearly $1<n_{max}<\infty$.

Let $g$ be an infinite order element of $G$ and $\w$ a cyclically reduced word in its conjugacy class. If $\w$ has
length at least $k$, the path $[\w^-,\w^+]$ is a $k$-local geodesic,
 and then a $(\l,\e)$-quasigeodesic. In particulary, for any $n\in \N$, one
has  $n|\w|\leq \l\,|\w^n|+\e$, and hence :
$$\frac{|\w|}{\l}-\frac{\e}{\l n}\leq \frac{|\w^n|}{n}$$
Making $n$ tend to $\infty$ one obtains $[[\w]]\geq |\w|/\l\geq k/\l$. Together with property \ref{1} one has :
$$[[g]]\geq \frac{k}{\l}$$

If $\w$ has length less than $k$, $\w$ lies in $\cal B$. Consider the integer $n_\w$ as defined above, then following
the same argument :
$$[[g^{n_\w}]]\geq \frac{k}{\l}\ \text{which implies}\ n_\w[[g]]\geq \frac{k}{\l} \text{
with property \ref{2}}$$ and finally for any infinite order element of $G$ :
$$[[g]]\geq K=\frac{k}{\l\,n_{max}}$$
which proves the assumption.\hfill$\square$

 \vskip 0.4cm


\end{document}